\documentclass[12pt]{article}
\usepackage{graphicx}
\usepackage{amsmath,amsthm,amssymb,enumerate}
\usepackage{euscript,mathrsfs}
\usepackage[left=2cm,right=2cm,top=3.5cm,bottom=3.5cm]{geometry}
\usepackage{color}

\usepackage{soul}

\catcode`\@=11 \@addtoreset{equation}{section}

\catcode`\@=12

\allowdisplaybreaks

\newtheorem{Theorem}{Theorem}[section]
\newtheorem{Proposition}[Theorem]{Proposition}
\newtheorem{Lemma}[Theorem]{Lemma}
\newtheorem{Corollary}[Theorem]{Corollary}

\theoremstyle{definition}
\newtheorem{Definition}[Theorem]{Definition}

\newtheorem{Remark}[Theorem]{Remark}

\newcommand{\bTheorem}[1]{
\begin{Theorem} \label{T#1} }
\newcommand{\eT}{\end{Theorem}}

\newcommand{\bProposition}[1]{
\begin{Proposition} \label{P#1}}
\newcommand{\eP}{\end{Proposition}}

\newcommand{\bLemma}[1]{
\begin{Lemma} \label{L#1} }
\newcommand{\eL}{\end{Lemma}}

\newcommand{\bCorollary}[1]{
\begin{Corollary} \label{C#1} }
\newcommand{\eC}{\end{Corollary}}

\newcommand{\bRemark}[1]{
\begin{Remark} \label{R#1} }
\newcommand{\eR}{\end{Remark}}

\newcommand{\bDefinition}[1]{
\begin{Definition} \label{D#1} }
\newcommand{\eD}{\end{Definition}}

\newcommand{\Del}{\Delta_x}

\newcommand{\bfphi}{\boldsymbol{\varphi}}

\newcommand{\bFormula}[1]{
\begin{equation} \label{#1}}
\newcommand{\eF}{\end{equation}}

\newcommand{\Ov}[1]{\overline{#1}}

\newcommand{\DC}{C^\infty_c}

\newcommand{\vr}{\varrho}

\newcommand{\vm}{\vc{m}}

\newcommand{\vc}[1]{{\bf #1}}

\newcommand{\Div}{{\rm div}_x}
\newcommand{\Grad}{\nabla_x}

\newcommand{\dx}{\,{\rm d} {x}}

\newcommand{\dt}{\,{\rm d} t }

\newcommand{\dxdt}{\dx  \dt}
\newcommand{\intO}[1]{\int_{\Omega} #1 \ \dx}

\newcommand{\ep}{\varepsilon}

\newcommand{\R}{\mathbb{R}}

\def\softd{{\leavevmode\setbox1=\hbox{d}%
          \hbox to 1.05\wd1{d\kern-0.4ex{\char039}\hss}}}
\definecolor{Cgrey}{rgb}{0.85,0.85,0.85}
\definecolor{Cblue}{rgb}{0.50,0.85,0.85}
\definecolor{Cred}{rgb}{1,0,0}
\definecolor{fancy}{rgb}{0.10,0.85,0.10}

\newcommand\Cbox[2]{%
    \newbox\contentbox%
    \newbox\bkgdbox%
    \setbox\contentbox\hbox to \hsize{%
        \vtop{
            \kern\columnsep
            \hbox to \hsize{%
                \kern\columnsep%
                \advance\hsize by -2\columnsep%
                \setlength{\textwidth}{\hsize}%
                \vbox{
                    \parskip=\baselineskip
                    \parindent=0bp
                    #2
                }%
                \kern\columnsep%
            }%
            \kern\columnsep%
        }%
    }%
    \setbox\bkgdbox\vbox{
        \color{#1}
        \hrule width  \wd\contentbox %
               height \ht\contentbox %
               depth  \dp\contentbox
        \color{black}
    }%
    \wd\bkgdbox=0bp%
    \vbox{\hbox to \hsize{\box\bkgdbox\box\contentbox}}%
    \vskip\baselineskip%
}

\date{}

\begin{document}


\title{On strong continuity of weak solutions to the compressible Euler system}

\author{Anna Abbatiello\thanks{The research of A.A. is supported by Einstein Foundation, Berlin.} \and Eduard Feireisl\thanks{The research of E.F. leading to these results has received funding from the
Czech Sciences Foundation (GA\v CR), Grant Agreement
18--05974S. The Institute of Mathematics of the Academy of Sciences of
the Czech Republic is supported by RVO:67985840. The stay of E.F. at TU Berlin is supported by Einstein Foundation, Berlin.} }


\maketitle

\bigskip

\centerline{Institute of Mathematics of the Academy of Sciences of the Czech Republic}
\centerline{\v Zitn\' a 25, CZ-115 67 Praha 1, Czech Republic}
\centerline{feireisl@math.cas.cz}
\centerline{and}

\centerline{Institute of Mathematics, Technische Universit\"{a}t Berlin,}
\centerline{Stra{\ss}e des 17. Juni 136, 10623 Berlin, Germany}
\centerline{anna.abbatiello@tu-berlin.de}

\begin{abstract}

Let $\mathcal{S} = \{ \tau_n \}_{n=1}^\infty \subset (0,T)$ be an arbitrary countable (dense) set. We show that for any given initial density and momentum, the compressible Euler system admits (infinitely many) admissible weak solutions that are not strongly continuous 
at each $\tau_n$, $n=1,2,\dots$. The proof is based on a refined version of the oscillatory lemma of De Lellis and Sz\' ekelyhidi 
with coefficients that may be discontinuous on a set of zero Lebesgue measure. 

\end{abstract}



\section{Introduction}
\label{I}

We consider the \emph{Euler system} describing the time evolution of the mass density $\vr = \vr(t,x)$ and the momentum $\vm = \vm(t,x)$ of 
a barotropic inviscid fluid: 
\begin{equation} \label{i1}
\begin{split}
\partial_t \vr + \Div \vm &= 0,\\
\partial_t \vm + \Div \left( \frac{\vm \otimes \vm}{\vr} \right) + \Grad p(\vr) &= 0, 
\end{split}
\end{equation}
where
$t \in (0,T)$, and $\Omega \subset \R^d$, $d=2,3$ is a bounded domain. The problem is supplemented by the impermeability condition 
\begin{equation} \label{i2}
\vm \cdot \vc{n}|_{\partial \Omega} = 0,
\end{equation}
and the initial conditions 
\begin{equation} \label{i3}
\vr(0, \cdot) = \vr_0, \ \vm(0, \cdot) = \vm_0.
\end{equation} 

As is well known, problem \eqref{i1}--\eqref{i3} is locally well posed for sufficiently regular initial data, however, the smooth 
solutions blow up in a finite time. The weak solution exists globally in time, however, the problem is essentially ill--posed 
even in the class of \emph{admissible weak solutions} satisfying the energy inequality 
\begin{equation} \label{i4}
\intO{ \left[ \frac{1}{2} \frac{ |\vm|^2}{\vr}  + P(\vr) \right] (t) } \leq 
\intO{ \left[ \frac{1}{2} \frac{ |\vm|^2}{\vr}  + P(\vr) \right] (s) } ,\ 
P'(\vr) \vr - P(\vr) = p(\vr)
\end{equation}
for a.a. $s$, including $s = 0$, and any $t$,\ $0 \leq s \leq t \leq T$. First examples of non--uniqueness were obtained in the 
seminal paper by DeLellis and Sz\' ekelyhidi \cite{DelSze3}, and later extended by Chiodaroli \cite{Chiod} and \cite{Fei2016}, 
Luo, Xie, and Xin \cite{LuXiXi} to a rather general class of initial data.

The key tool for using the convex integration machinery of \cite{DelSze3}, developed originally for the incompressible fluids, is 
a suitable adaptation of the so-called \emph{Oscillatory Lemma}, proved originally in \cite{DelSze3} and extended to ``variable coefficients'' in \cite{Chiod}. Probably the most general version including ``non--local coefficients'' can be found 
in \cite{Fei2016}. The limitation of this approach is due to the fact that certain quantities, in particular the initial density 
and the desired energy profile, must enjoy some degree of smoothness to transform the problem to its basic form handled 
in \cite{DelSze3}. The largest possible class used so far is that of \emph{piecewise continuous} functions, cf. \cite{Fei2016}, 
\cite{LuXiXi}. 

A closer inspection of the problem reveals apparent similarity between the regularity properties required 
for the coefficients in Oscillatory Lemma and their \emph{integrability} in the Riemann sense. Our goal is to extend validity 
of Oscillatory Lemma to the case of Riemann integrable coefficients, specifically belonging to the class:
\[
\mathcal{R}(Q) \equiv \left\{ \ v: Q \to \R \ \Big|\ 
{\rm meas} \left\{ y \in Q \ \Big| \ v \ \mbox{is not continuous at} \ y     \right\} = 0   \right\}
\]     
where the symbol ``meas" stands for the Lebesgue measure.
Such an extension allows us to show the existence of weak solutions to the Euler system with a given total energy profile belonging 
to $\mathcal{R}$. In particular, as the weak solutions $[\vr, \vm]: t \mapsto L^1(\Omega) \times L^1(\Omega; \R^d)$ 
are strongly continuous at a time $t$ if and only if the total energy is continuous at $t$, we obtain the existence of 
an admissible weak solution 
that is not strongly continuous 
at an arbitrary given countable dense set of times $\mathcal{S} = \{ \tau_n \}_{n=1}^\infty \subset (0,T)$. 

The paper is organized as follows. In Section \ref{M}, we collect the preliminary material and state our main results. In Section 
\ref{O}, we show a version of Oscillatory Lemma with coefficients belonging to $\mathcal{R}$. Applications, including the proofs of the 
main results, are discussed in Section \ref{A}.

\section{Preliminaries, main results}
\label{M}

We say that the functions 
\[
\vr \in C_{{\rm weak}}([0,T]; L^2(\Omega)),\ \vm \in C_{{\rm weak}}([0,T]; L^2(\Omega; \R^d))
\]
represent \emph{weak solution} to the Euler problem \eqref{i1}--\ref{i3} if:
\begin{itemize}
\item $\vr \geq 0$, $p(\vr) \in L^1((0,T) \times \Omega)$;
\item
the equation of continuity
\begin{equation} \label{M1}
\int_0^T \intO{ \left[ \vr \partial_t \varphi + \vm \cdot \Grad \varphi \right]   } \dt = - \intO{ \vr_0 \varphi (0, \cdot) }
\end{equation}
holds
for any $\varphi \in C^1_{{\rm loc}}([0,T) \times \Ov{\Omega})$;
\item
the momentum equation
\begin{equation} \label{M2}
\int_0^T \intO{ \left[ \vm \cdot \partial_t \bfphi + \left( \frac{\vm \otimes \vm}{\vr} \right) : \Grad \bfphi  
+ p(\vr) \Div \bfphi \right]   } \dt = - \intO{ \vm_0 \bfphi (0, \cdot) }
\end{equation}
holds
for any $\bfphi \in C^1_{{\rm loc}}([0,T) \times \Ov{\Omega}; \R^d)$, $\bfphi \cdot \vc{n} |_{\partial \Omega} = 0$.
\end{itemize}

A weak solution $[\vr, \vm]$ is \emph{admissible} if it satisfies the energy inequality \eqref{i4} for any $t \in (0,T)$ and a.a.  
$s \in (0,T)$, $0 \leq s < t$.  

\subsection{Main results, solutions with arbitrary energy profile}

We are ready to state our first result. 

\begin{Theorem} \label{MT1}
Let $\Omega \subset \R^d$, $d = 2,3$, be a bounded domain with $C^2$ boundary. Let the initial data $\vr_0$, $\vm_0$ be given,
\[
0 < \underline{\vr} \leq \vr_0(x) \leq \Ov{\vr} \ \mbox{for all}\ x \in \Ov{\Omega}, \ \vr_0 \in \mathcal{R}(\Ov{\Omega}),
\]
\[
\vm_0 \in \mathcal{R}(\Ov{\Omega}; R^d),\ \Div \vm_0 \in \mathcal{R}(\Ov{\Omega}), \ \vm_0 \cdot \vc{n}|_{\partial \Omega} = 0.
\]
Let $E(t)$ be an arbitrary function satisfying 
\[
0 \leq E(t) \leq \Ov{E} \ \mbox{for all}\ t \in [0,T], \ E \in \mathcal{R}[0,T].
\]

Then there exist $E_0 \geq 0$ such that the Euler system 
\eqref{M1}, \eqref{M2} admits infinitely many solutions $[\vr, \vm]$ in $(0,T) \times \Omega$ satisfying 
\[
\frac{1}{2}\underline{\vr} \leq \vr(t,x) \leq 2 \Ov{\vr}\ \mbox{for all } (t,x) \in (0,T) \times \Omega ,\ \vm \in L^\infty((0,T) \times \Omega; \R^d)),
\]
\[
\intO{ \left[ \frac{1}{2} \frac{|\vm|^2}{\vr} + P(\vr) \right](\tau, \cdot) } = E_0 + E(\tau) \ \mbox{for a.a.}\ \tau \in (0,T).  
\]
 
\end{Theorem}

\begin{Remark} \label{RR1}

It will be clear from the proof given below that 
the density profile can be taken $\vr = \vr_0(x)$ as soon as $\Div \vm_0 = 0$. In such a case, we may consider 
$\vr_0 \equiv 1$ obtaining the same conclusion for the incompressible Euler system. Moreover, the result holds for 
any bounded domain, no smoothness of the boundary is necessary. 

\end{Remark}

Solutions satisfying strict energy inequality cannot be regular, cf. e.g. Constantin, E, and Titi \cite{ConETit} 
or \cite{FeGwGwWi}.  
Similarly to other ``wild'' solutions 
produced by the method of convex integration, the solutions may experience the initial energy jump, meaning 
the energy inequality \eqref{i4} may not hold for $s=0$. 
However, as there is definitely 
a sequence of times $\tau_n \searrow 0$ for which 
\[
\intO{ \left[ \frac{1}{2} \frac{|\vm|^2}{\vr} + P(\vr) \right](\tau_n, \cdot) } = E_0 + E(\tau_n).
\] 
One could also deduce the existence of infinitely many solutions with the energy continuous at the initial time\color[rgb]{0,0,0}, performing the procedure 
described e.g. in \cite{DelSze3}. We leave the details to the interested reader.

\subsection{Strong continuity in time}

We say that a weak solution $[\vr, \vm]$ of the Euler system is \emph{strongly continuous} at a time $\tau \in (0,T)$ if 
\[
\vr(t, \cdot) \to \vr(\tau, \cdot)\ \mbox{in}\ L^1(\Omega),\ 
\vm(t, \cdot) \to \vm(\tau, \cdot)\ \mbox{in}\ L^1(\Omega; \R^d) \ \mbox{for}\ t \to \tau.
\]

\begin{Theorem} \label{MT2}
Let $\Omega \subset \R^d$, $d = 2,3$, be a bounded domain with $C^2$ boundary. 
Let the initial data $\vr_0$, $\vm_0$ be given,
\[
0 < \underline{\vr} \leq \vr_0(x) \leq \Ov{\vr} \ \mbox{for all}\ x \in \Ov{\Omega}, \ \vr_0 \in \mathcal{R}(\Ov{\Omega}),
\]
\[
\vm_0 \in \mathcal{R}(\Ov{\Omega}; R^d),\ \Div \vm_0 \in \mathcal{R}(\Ov{\Omega}), \ \vm_0 \cdot \vc{n}|_{\partial \Omega} = 0.
\]
Let $\mathcal{S} = \{ \tau_n \}_{n=1}^\infty \subset (0,T)$ be an arbitrary (countable) set of times. 

Then the Euler system admits infinitely many admissible weak solutions  
that are not strongly continuous at any $\tau_n$, $n=1,2,\dots$  

\end{Theorem}

Here again \emph{admissible} means  the total energy is equal to a non--increasing function
for a.a. time. In particular, the solutions need not be strongly continuous at $t = 0$.

\section{Oscillatory lemma}
\label{O}

The proof of our main results depends on a generalized version of Oscillatory Lemma of De Lellis and Sz\' ekelyhidi \cite{DelSze3}. 
Our starting point is its most elementary version showed in \cite[Proposition 3]{DelSze3}:

\begin{Lemma}[Oscillatory Lemma, basic form] \label{OL1}

Let $Q = (0,1) \times (0,1)^d$, $d=2,3$. Suppose that $\vc{v} \in \R^d$, $\mathbb{U} \in \R^{d \times d}_{0,{\rm sym}}$, 
$e \leq \Ov{e}$ are given 
constant quantities satisfying\footnote{$\R^{d \times d}_{0,{\rm sym}}$ denotes the space of real symmetric matrices with zero trace, while 
$\lambda_{\rm max}[\cdot]$ is the maximum eigenvalue.}  
\[
\frac{d}{2} \lambda_{\rm max} \left[ \vc{v} \otimes \vc{v} - \mathbb{U} \right] < e.
\] 

Then there is a constant $c = c(d, \Ov{e})$ and sequences of vector functions $\{ \vc{w}_n \}_{n=1}^\infty$, 
$\{ \mathbb{V}_n \}_{n=1}^\infty$, 
\[
\vc{w}_n \in \DC(Q; \R^d),\ \mathbb{V}_n \in \DC(Q; \R^{d \times d}_{0, {\rm sym}}) 
\]
satisfying 
\[ \partial_t \vc{w}_n + \Div \mathbb{V}_n = 0,\ \Div \vc{w}_n = 0 \mbox{ in } Q,\]
\[
\frac{d}{2} \lambda_{\rm max} \left[ (\vc{v} + \vc{w}_n) \otimes (\vc{v} + \vc{w}_n) - (\mathbb{U} + 
\mathbb{V}_n ) \right] < e \ \mbox{in}\ Q \ \mbox{for all}\ n = 1,2, \dots,
\]
\[
\vc{w}_n \to 0 \ \mbox{in}\ C_{\rm weak}([0,1]; L^2((0,1)^d; \R^d))\ \mbox{as}\ n \to \infty,
\]
\[
\liminf_{n \to \infty} \int_Q |\vc{w}_n|^2 \dxdt \geq c(d, \Ov{e}) \int_Q \left( e - \frac{1}{2} |\vc{v}|^2 \right)^2 \dxdt.
\]

\end{Lemma} 

\subsection{Extension by scaling}

We say that $Q \subset [0,T] \times \R^d$ is a \emph{block}, if 
\[
Q = (t_1, t_2) \times \Pi_{i=1}^d (a_i, b_i),\ t_1 < t_2,\ a_i < b_i, \ i=1,\dots,d.
\]
The following can be easily deduced from Lemma \ref{OL1} by a scaling argument, see e.g. Chiodaroli \cite[Section 6, formula (6.9)]{Chiod}.

\begin{Lemma}[Oscillatory Lemma, scaled form] \label{OL2}

Let 
\[
Q = (t_1, t_2) \times \Pi_{i=1}^d (a_i, b_i),\ t_1 < t_2,\ a_i < b_i, \ i=1,\dots,d,
\]
be a block.
Suppose that $\vc{v} \in \R^d$, $\mathbb{U} \in \R^{d \times d}_{0,{\rm sym}}$, 
$e \leq \Ov{e}$, and $r > 0$ are given 
constant quantities satisfying
\[
\frac{d}{2} \lambda_{\rm max} \left[ \frac{\vc{v} \otimes \vc{v}}{r} - \mathbb{U} \right] < e.
\] 

Then there is a constant $c = c(d, \Ov{e})$ and sequences of vector functions $\{ \vc{w}_n \}_{n=1}^\infty$, 
$\{ \mathbb{V}_n \}_{n=1}^\infty$, 
\[
\vc{w}_n \in \DC(Q; \R^d),\ \mathbb{V}_n \in \DC(Q; \R^{d \times d}_{0, {\rm sym}}) 
\]
satisfying 
\[ \partial_t \vc{w}_n + \Div \mathbb{V}_n = 0,\ \Div \vc{w}_n = 0 \mbox{ in } Q,\]
\[
\frac{d}{2} \lambda_{\rm max} \left[ \frac{(\vc{v} + \vc{w}_n) \otimes (\vc{v} + \vc{w}_n)}{r} - (\mathbb{U} + 
\mathbb{V}_n ) \right] < e\ \mbox{in}\ Q\ \mbox{for all}\ n = 1,2, \dots,
\]
\[
\vc{w}_n \to 0 \ \mbox{in}\ C_{\rm weak}([t_1, t_2]; L^2(\Pi_{i=1}^d (a_i, b_i); \R^d))\ \mbox{as}\ n \to \infty,
\]
\[
\liminf_{n \to \infty} \int_Q \frac{|\vc{w}_n|^2}{r} \dxdt \geq c(d, \Ov{e}) \int_Q \left( e - \frac{1}{2} \frac{|\vc{v}|^2}{r} \right)^2 \dxdt.
\]

\end{Lemma}

\subsection{Oscillatory Lemma for Riemann integrable coefficients}

Our main goal is to show the following extension of Oscillatory Lemma. 

\begin{Lemma}[Oscillatory Lemma, general coefficients] \label{OL3}

Let 
\[
Q = (t_1, t_2) \times \Pi_{i=1}^d (a_i, b_i),\ t_1 < t_2,\ a_i < b_i, \ i=1,\dots,d,
\]
be a block.
Suppose that 
\[
\vc{v} \in \mathcal{R}(\Ov{Q}; \R^d),\ \mathbb{U} \in \mathcal{R}(\Ov{Q}; \R^{d \times d}_{0,{\rm sym}}), 
\ e \in \mathcal{R}(\Ov{Q}),\ r \in \mathcal{R}(\Ov{Q})
\]
be given such that
\begin{equation} \label{HYP1}
0 < \underline{r} \leq r(t,x) \leq \Ov{r},\ e(t,x) \leq \Ov{e} \ \mbox{for all}\ (t,x) \in \Ov{Q}, 
\end{equation}
\begin{equation} \label{HYP2}
\frac{d}{2} \sup_{\Ov{Q}}
\lambda_{\rm max} \left[ \frac{\vc{v} \otimes \vc{v}}{r} - \mathbb{U} \right] < \inf_{\Ov{Q}} e.
\end{equation}

Then there is a constant $c = c(d, \Ov{e})$ and sequences of vector functions $\{ \vc{w}_n \}_{n=1}^\infty$, 
$\{ \mathbb{V}_n \}_{n=1}^\infty$, 
\[
\vc{w}_n \in \DC(Q; \R^d),\ \mathbb{V}_n \in \DC(Q; \R^{d \times d}_{0, {\rm sym}}) 
\]
satisfying 
\[ \partial_t \vc{w}_n + \Div \mathbb{V}_n = 0,\ \Div \vc{w}_n = 0 \mbox{ in } Q,\]
\begin{equation} \label{CONCL1}
\frac{d}{2} \sup_{\Ov{Q}} \lambda_{\rm max} \left[ \frac{(\vc{v} + \vc{w}_n) \otimes (\vc{v} + \vc{w}_n)}{r} - (\mathbb{U} + 
\mathbb{V}_n ) \right] < \inf_{\Ov{Q}} {e}\ \mbox{for all}\ n = 1,2, \dots,
\end{equation}
\[
\vc{w}_n \to 0 \ \mbox{in}\ C_{\rm weak}([t_1, t_2]; L^2(\Pi_{i=1}^d (a_i, b_i); \R^d))\ \mbox{as}\ n \to \infty,
\]
\[
\liminf_{n \to \infty} \int_Q \frac{|\vc{w}_n|^2}{r} \dxdt \geq c(d, \Ov{e}) \int_Q \left( e - \frac{1}{2} \frac{|\vc{v}|^2}{r} \right)^2 \dxdt.
\]

\end{Lemma}

The remaining part of this section will be devoted to the proof of Lemma \ref{OL3}.

\subsubsection{Basic properties of Riemann integrable functions} 

The leading idea is to approximate the coefficients $\vc{v}$, $\mathbb{U}$, $e$, and $r$ by piecewise constant functions 
and use Lemma \ref{OL2}. The following is standard and may be found e.g. in the textbook by Zorich \cite[Chapter 11]{Zor}.
 
For a real valued function $v : \Ov{Q} \to \R$ we introduce:
\[
{\rm osc} [v](t,x) = \lim_{s \searrow 0} \left[ \sup_{ B((t,x), s) \cap \Ov{Q} } v -  \inf_{ B((t,x), s) \cap \Ov{Q} } v \right],
\]
where $B((t,x), s)$ denotes the ball of radius $s$ centered at $(t,x)$. 
It holds:

\begin{itemize}

\item

\begin{equation} \label{S2}
A_\eta = \left\{ (t,x) \in \Ov{Q}\ \Big|\ {\rm osc} [v] (t,x) \geq \eta \right\} \ \mbox{is closed}
\end{equation} 

\item for any $v \in \mathcal{R}(\Ov{Q})$ and $\eta > 0$, the set $A_\eta$ is of zero content, meaning for any $\delta > 0$, there exists a 
\emph{finite} number of (open) boxes $Q_i$ such that 
\[
A_{\eta} \subset \cup_i Q_i,\ \sum_i |Q_i| < \delta.
\]

\end{itemize}

\subsubsection{Continuity of eigenvalues} 

We recall the algebraic inequalities (see e.g. \cite{DelSze3}) 
\begin{equation} \label{O1}
\frac{1}{2} \frac{|\vc{v}|^2} {r} \leq \frac{d}{2} \lambda_{\rm max} \left[ \frac{\vc{v}  \otimes \vc{v} }{r} - \mathbb{U}  \right],
\  \|\mathbb{U}\|_\infty\leq \frac{d}{2} \lambda_{\rm max} \left[ \frac{\vc{v}  \otimes \vc{v} }{r} - \mathbb{U}  \right]
\end{equation}
for any $\vc{v} \in \R^d,$  $r > 0,$  $\mathbb{U} \in \R^{d \times d}_{0,{\rm sym}}$, where $\|\mathbb{U}\|_\infty$ denotes the operator norm of the matrix.\\
Consider the set
\[
K = \left\{ r \in (0, \infty), \vc{v} \in \R^d, \mathbb{U} \in \R^{d \times d}_{0,{\rm sym}} \ 
\Big|\ \underline{r} \leq r \leq \Ov{r},\ \frac{d}{2} \lambda_{\rm max} \left[ \frac{\vc{v}  \otimes \vc{v} }{r} - \mathbb{U}  \right] 
\leq \Ov{e}  \right\}.
\]
In view of \eqref{O1}, $K$ is a compact subset of $(0, \infty) \times \R^d \times \R^{d \times d}_{0,{\rm sym}}$. Moreover, as shown in 
\cite{DelSze3}, the function 
\[
[\vc{w}, \mathbb{U}] \mapsto \frac{d}{2} \lambda_{\rm max} \left[ {\vc{w}  \otimes \vc{w} } - \mathbb{U}  \right]
\]
is convex. As convex functions are Lipschitz continuous on compact subsets of their domain, we deduce there is a constant $L$ such that 
\begin{equation} \label{O2}
\begin{split}
\left| \lambda_{\rm max} \left[ \frac{\vc{v}_1  \otimes \vc{v}_1 }{r_1} - \mathbb{U}_1  \right] 
- \lambda_{\rm max} \left[ \frac{\vc{v}_2  \otimes \vc{v}_2 }{r_2} - \mathbb{U}_2 \right] \right| &\leq L
\Big( |r_1 - r_2| + |\vc{v}_1 - \vc{v}_2| + |\mathbb{U}_1 - \mathbb{U}_2 | \Big)\\
&\mbox{for any}\ (r_i, \vc{v}_i, \mathbb{U}_i) \in K,\ i=1,2.
\end{split}
\end{equation}

\subsubsection{Domain decomposition}

Suppose $\vc{v}$, $r$, $\mathbb{U}$, $r$, and $e$ satisfy \eqref{HYP1}, \eqref{HYP2}. It follows from \eqref{HYP2} 
that there exists $\ep_0 > 0$ such that 
\[ 
\frac{d}{2} \lambda_{\rm max} \left\{ \frac{\vc{v} \otimes \vc{v}}{r} - \mathbb{U} \right\} < e - \ep_0 \leq \Ov{e} \ \mbox{in}\ \Ov{Q}.
\]
In particular $(r, \vc{v}, \mathbb{U})(t,x) \in K$ for any $(t,x) \in \Ov{Q}$.
Thus for any $0 < \ep \leq \ep_0$  
\begin{equation} \label{S3}
\frac{d}{2} \lambda_{\rm max} \left\{ \frac{\vc{v} \otimes \vc{v}}{r} - \mathbb{U} \right\} < e - \ep \leq 
\Ov{e}\  \mbox{in}\ \Ov{Q}.
\end{equation}

For $\eta > 0$ consider the set 
\[
A_\eta \equiv 
A_{\eta}[\vc{v}] \cup A_{\eta} [r] \cup A_{\eta}[\mathbb{U}] \cup A_{\eta} [e],
\]  
cf. \eqref{S2}.
In accordance with our hypotheses, this is a set of zero content, meaning 
there is a finite number of (open) boxes $Q^s_i(\eta)$ such that
\[
A_{\eta} \subset \cup_i Q^s_i(\eta),\  
\sum_i |Q^s_i(\eta)| < \ep \ \mbox{for any given}\ \eta > 0.
\]

The complement $\Ov{Q} \setminus \cup_i Q^s_i(\eta)$ is compact. Moreover, each point $y \in \Ov{Q} \setminus \cup_i Q^s_i$ 
has an open neighborhood $U(y)$ such that 
\begin{equation} \label{S5}
|r(y_1) - r(y_2)| < 2\eta,\ |\vc{v}(y_1) - \vc{v}(y_2)| < 2\eta,\ 
|\mathbb{U}(y_1) - \mathbb{U}(y_2)| < 2\eta, \ |e(y_1) - e(y_2)| <2 \eta   
\end{equation}
whenever $y_1, y_2 \in U(y) \cap \Ov{Q}$.

As the set $\Ov{Q} \setminus \cup_i Q^s_i(\eta)$ is compact and there is a finite number of 
$Q^s_i$, we may infer that for any 
given $\ep > 0$, $\eta > 0$, there exists a decomposition of $\Ov{Q}$ into a finite number of blocks: 
\begin{equation} \label{S6}
\begin{split}
\Ov{Q} = (\cup_{i} \Ov{Q}^s_i(\eta)) &\cup (\cup_j \Ov{Q}^r_j(\eta) ),\ Q^r_j \cap Q^r_k = \emptyset 
\ \mbox{for}\ j \ne k,\\
\mbox{such that} &\\
\sum_{i} |Q^s_i(\eta)| &< \ep\\
\mbox{and}&\\
|r(y_1) - r(y_2)| &<2 \eta,\ |\vc{v}(y_1) - \vc{v}(y_2)| <2 \eta,\ 
|\mathbb{U}(y_1) - \mathbb{U}(y_2)| < 2\eta,\ |e(y_1) - e(y_2)| <2 \eta\\
&\mbox{for any}\ y_1, y_2 \in \Ov{Q}^r_j,\ j = 1,2,\dots
\end{split}
\end{equation}

\subsubsection{Localization}

Given $0 < \ep \leq \ep_0$, $\eta > 0$, consider the decomposition of $\Ov{Q}$ given by \eqref{S6}. Choosing $y_j \in Q^r_j$ we fix 
\[
\widetilde{r} = r(y_j), \ \widetilde{\vc{v}} = \vc{v}(y_j),\  \widetilde{\mathbb{U}} = \mathbb{U} (y_j), \ \mbox{and}\ 
\widetilde{e} = e(y_j). 
\]
Applying the constant coefficient version of Oscillatory Lemma (Lemma \ref{OL2}) on each $Q^r_j$
we get a sequence of functions $\vc{w}^j_n$, $\mathbb{V}^j_n$, smooth and compactly supported in $Q^r_j \equiv (s_1, s_2)\times O^r_j$, such that
\begin{equation} \label{S8}
\partial_t \vc{w}^j_n + \Div \mathbb{V}^j_n = 0,\ \Div \vc{w}^j_n = 0 \mbox{ in } Q^r_j,
\end{equation}
\begin{equation} \label{S9}
\vc{w}^j_n \to 0 \ \mbox{in}\ C_{\rm weak}([s_1,s_2]; L^2(O^r_j; \R^d)),
\end{equation}
\begin{equation} \label{S10}
\frac{d}{2} \lambda_{\rm max} \left\{ \frac{(\widetilde{\vc{v}} + \vc{w}^j_n) \otimes (\widetilde{\vc{v}} + \vc{w}^j_n)}{
\widetilde{r}} - 
(\widetilde{\mathbb{U}} + \mathbb{V}^j_n ) \right\} + \ep < \widetilde{e} \mbox{ in } Q^r_j,
\end{equation}
and
\begin{equation} \label{S11}
\liminf_{n \to \infty} \int_{Q^r_j} \frac{|\vc{w}^j_n|^2}{\widetilde{r}} \ \dxdt  \geq c(d, \Ov{e}) \int_{Q^r_j} \left( \widetilde{e} - \frac{1}{2} \frac{|\widetilde{\vc{v}}|^2}{\widetilde{r}} \right)^2 \dxdt.
\end{equation}

In view of the Lipschitz continuity of the eigenvalues established in \eqref{O2}, and in accordance with \eqref{S6}, we may choose 
$\eta = \eta(\ep)$ small enough so that 
\begin{equation} \label{S12}
\frac{d}{2} \lambda_{\rm max} \left\{ \frac{({\vc{v}} + \vc{w}^j_n) \otimes ({\vc{v}} + \vc{w}^j_n)}{
{r}} - 
({\mathbb{U}} + \mathbb{V}^j_n ) \right\} + \frac{\ep}{2} < {e} \ \mbox{in}\ Q^r_j.
\end{equation}
By the same token, we get 
\begin{equation} \label{S13}
\liminf_{n \to \infty} \int_{Q^r_j} \frac{|\vc{w}^j_n|^2}{{r}} \ \dxdt  \geq c(d, \Ov{e}) \int_{Q^r_j} \left( {e} - \frac{1}{2} \frac{|{\vc{v}}|^2}{{r}} \right)^2 \ \dxdt - \ep| Q^r_j |.
\end{equation}

Finally, setting $\vc{w}^i_n = \mathbb{V}^i_n = 0$ on $Q^s_i$ and summing up over all boxes, we obtain
sequences defined on $\Ov{Q}$ satisfying 
\begin{equation} \label{S14}
\begin{split}
\liminf_{n \to \infty} \int_{Q} \frac{|\vc{w}_n|^2}{{r}} \ \dxdt  &\geq c(d, \Ov{e}) 
\sum_{j} \int_{Q^r_j} \left( {e} - \frac{1}{2} \frac{|{\vc{v}}|^2}{{r}} \right)^2 \ \dxdt - \ep |Q|\\
&\geq c(d, \Ov{e}) 
\int_{Q} \left( {e} - \frac{1}{2} \frac{|{\vc{v}}|^2}{{r}} \right)^2 \ \dxdt - \ep \left(|Q| + \Ov{e}^2 \right),
\end{split}
\end{equation}
and 
\begin{equation}\label{ann}
\vc{w}_n \to 0 \ \mbox{in}\ C_{\rm weak}([t_1, t_2]; L^2(\Pi_{i=1}^d (a_i, b_i); \R^d)).
\end{equation}
As pointed out, the oscillatory perturbations can be constructed for any $0 < \ep < \ep_0$.

\subsubsection{Diagonalization argument}

To complete the proof of Lemma \ref{OL3}, it remains to get rid of the $\ep-$dependent term in \eqref{S14}.
This can be achieved by a simple diagonalization argument. 
By the previous subsection, for any $\ep>0$ there exists $\{\vc{w}_n^\ep\}_{n\in \mathbb{N}}$ such that \eqref{S14} and \eqref{ann} hold. 
Combining \eqref{S14} and a basic property of the liminf, we get that  there exists $n_{0,\ep}$ such that for all $n\geq n_{0,\ep}$ it holds
\begin{equation}\label{liminf0}
 \int_{Q} \frac{|\vc{w}_n^\ep|^2}{{r}} +\ep (|Q| + \Ov{e}^2)\ \dxdt  \geq c(d, \Ov{e}) \int_{Q} \left( {e} - \frac{1}{2} \frac{|{\vc{v}}|^2}{{r}} \right)^2 \ \dxdt.
\end{equation}
In addition we can fix $n_{0,\ep}$ in such a way that 
\begin{equation}\label{lim} d(\vc{w}_n^\ep, 0) < \ep \mbox{ for all } n\geq n_{0,\ep}\end{equation}
where $d(\cdot,\cdot)$ is the metric defined as
$$d(\cdot,\cdot)= \sup_{t\in [0,T]} d_B (\cdot,\cdot)$$
and $d_B (\cdot,\cdot)$ is the metric  induced by the weak topology on bounded sets of the Hilbert space $L^2(\Pi_{i=1}^d (a_i, b_i); \R^d)$.
For any $k\in \mathbb{N}$, let us choose $\ep=\frac{\ep_0}{k}$ then there exists a sequence $ \{\vc{w}_n^{\frac{1}{k}}\}_{n\in \mathbb{N}}$, which fulfills \eqref{liminf0} and \eqref{lim} definitely. We do not relabel such subsequence. Thus we get an infinite matrix
 $$\begin{pmatrix}
 \vc{w}_1^1 & \vc{w}_2^1 & \cdots & \vc{w}_k^1 &\cdots & \vc{w}_n^1&\cdots \\
\vc{w}_1^{1/2} & \vc{w}_2^{1/2} & \cdots & \vc{w}_k^{1/2} & \cdots &  \vc{w}_n^{1/2} &\cdots\\
  \vdots  & \vdots  & \ddots & \vdots &\ddots &\vdots&\cdots\\
\vc{w}_1^{1/k} & \vc{w}_2^{1/k} & \cdots & \vc{w}_k^{1/k} &\cdots &  \vc{w}_n^{1/k}&\cdots\\
  \vdots  & \vdots  & \ddots & \vdots& \ddots &\vdots&\cdots
 \end{pmatrix}$$
 where the $k$-th row corresponds to a sequence fulfilling \eqref{liminf0} and \eqref{lim} with $\ep=\frac{\ep_0}{k}$. Consider the sequence $\{\vc{w}_k^{1/k}\}_k$, which corresponds to the diagonal of the matrix above, it enjoys 
 \begin{equation}\label{liminf1}
 \int_{Q} \frac{|\vc{w}_k^{1/k}|^2}{{r}} +\frac{1}{k} (|Q| + \Ov{e}^2)\ \dxdt  \geq c(d, \Ov{e}) \int_{Q} \left( {e} - \frac{1}{2} \frac{|{\vc{v}}|^2}{{r}} \right)^2 \ \dxdt
\end{equation}
and 
\begin{equation}
d(\vc{w}_k^{1/k}, 0)<\frac{\ep_0}{k}.
\end{equation}
Taking respectively the liminf and the limit as $k\rightarrow+\infty$, we conclude that 
 \begin{equation}\label{liminf2}
\liminf_{k\to+\infty} \int_{Q} \frac{|\vc{w}_k^{1/k}|^2}{{r}}  \dxdt  \geq c(d, \Ov{e}) \int_{Q} \left( {e} - \frac{1}{2} \frac{|{\vc{v}}|^2}{{r}} \right)^2 \ \dxdt
\end{equation}
and 
\begin{equation}
\vc{w}_k^{1/k} \to 0 \mbox{ in } C_{\rm weak}([t_1, t_2]; L^2(\Pi_{i=1}^d (a_i, b_i); \R^d)).
\end{equation}

\begin{Remark} \label{OR1}

The conclusion of Lemma \ref{OL3} holds if $Q$ is a bounded open set. Indeed $Q$ can be covered by a countable number of blocks on each of which we may apply the previous arguments. 
The relevant result is provided by Whitney decomposition lemma (Stein \cite{STEIN1}),
see \cite[Section 4.4]{DoFeMa} for details.

\end{Remark}

\section{Applications}
\label{A}

Our ultimate goal is to apply the general version of Oscillatory Lemma to show existence of weak solutions to the compressible Euler system with given energy. 

\subsection{Rewriting the Euler system as an abstract problem}

Following \cite{Fei2016}, we write the initial momentum in the form of its Helmholtz decomposition, 
\[
\vm_0 = \vc{v}_0 + \Grad \Phi_0,
\] 
where 
\[
\Del \Phi_0 = \Div \vm_0 \ \mbox{in}\ \Omega, \ \Grad \Phi_0 \cdot \vc{n}|_{\partial \Omega} = 0.
\]
As the boundary $\partial \Omega$ is of class $C^2$, the standard elliptic estimates imply
$\Grad \Phi_0 \in W^{1,p}(\Omega; R^d)$, in particular $\Grad \Phi_0 \in C(\Ov{\Omega}; R^d)$, see e.g. Agmon, Douglis, and Nirenberg  
\cite{ADN}.

Next, we fix the density profile 
\[
\vr(t,x) = \vr_0 + h(t) \Del \Phi_0, \ h \in C^\infty[0, \infty), \ h(0) = 0, \ h'(0) = -1.
\]
We look for solutions in the form 
\[
\vc{m} = \vc{v} - h'(t) \Grad \Phi_0, \ \Div \vc{v} = 0, \ \vc{v} \cdot \vc{n}|_{\partial \Omega} = 0.
\]
Seeing that 
\[
\partial_t \vr = h'(t) \Del \Phi_0 = - \Div \vc{m}, 
\]
we can adjust $h$ in such a way that 
\[
0 < \frac{1}{2} \underline{\vr} \leq \vr(t,x) \leq 2 \Ov{\vr} \ \mbox{for all}\ t \geq 0,\ x \in \Ov{\Omega}
\]
provided the initial density is uniformly bounded below and above. In addition, for $\vr_0 \in \mathcal{R}(\Ov{\Omega})$, 
we have 
\[
\vr \in \mathcal{R}([0,T] \times \Ov{\Omega}).
\]

Accordingly, we look for a vector field $\vc{v}$ solving the following problem:
\begin{equation} \label{S15}
\begin{split}
&\Div \vc{v} = 0, \ \vc{v} \cdot \vc{n}|_{\partial \Omega} = 0,\ \vc{v}(0, \cdot) = \vc{v}_0, 
\\
&\partial_t \vc{v} - h''(t) \Grad \Phi_0  + \Div \left( \frac{(\vc{v} - h'(t) \Grad \Phi_0) \otimes (\vc{v}  - 
h'(t) \Grad \Phi_0) }{\vr} - \frac{1}{d} \frac{|\vc{v}  - h'(t) \Grad \Phi_0 |^2}{\vr} \mathbb{I} \right) = 0,
\end{split}
\end{equation}
with prescribed kinetic energy
\begin{equation} \label{S16}
\frac{1}{2} \frac{|\vc{v} - h'(t)\Grad \Phi_0|^2}{\vr} = \Lambda(t) - \frac{d}{2} p(\vr) + \frac{d}{2} h''(t) \Phi_0
\end{equation}
where $\Lambda = \Lambda(t)$ is a spatially homogeneous function to be chosen below.

Obviously (cf. Chiodaroli \cite{Chiod} and \cite{Fei2016}), any \emph{weak} solution $\vc{v}$ of \eqref{S15}, \eqref{S16} gives rise to a weak solution $[\vr, \vc{m} = \vc{v} - h'(t) \Grad \Phi_0]$ of the Euler system \eqref{M1}, \eqref{M2}, with the total energy 
\begin{equation} \label{S17}
\intO{ \left[ \frac{1}{2} \frac{|\vm  |^2}{\vr } + P(\vr) \right](\tau, \cdot)} = 
\Lambda (\tau) |\Omega| + \intO{ \left[ P(\vr) - \frac{d}{2} p(\vr) + \frac{d}{2} h''(\tau) \Phi_0 \right] }\ \mbox{for a.a.}\ \tau \in (0,T).
\end{equation}
 
Evoking the notation of Theorem \ref{MT1}, we set
\[
\Lambda(\tau) = \frac{E(\tau)}{|\Omega|} + \Lambda_0 (\tau) ,\ E_0 = 
\Lambda_0(\tau) |\Omega| + \intO{ \left[ P(\vr) - \frac{d}{2} p(\vr) + h''(\tau) \Phi_0 \right] }.
\]
Thus the proof of Theorem \ref{MT1} consists in showing that for given $\vr_0$ and $E$, there exists $E_0$ large enough 
so that the problem \eqref{S15}, \eqref{S16} admits (infinitely many) weak solutions.  

\subsection{Subsolutions}

We start by fixing the energy profile
\[
e = e(t,x) = \frac{E(t)}{|\Omega|} + \Lambda_0(t) - \frac{d}{2} p(\vr) + \frac{d}{2} 
h''(t) \Phi_0,\ e \in \mathcal{R}([0,T] \times \Ov{\Omega}).
\]
Similarly to \cite{DelSze3}, we introduce the space of \emph{subsolutions},
\[
\begin{split}
X_0 = &\left\{ (\vc{v} - \vc{v}_0) \in C^1([0,T] \times \Ov{\Omega}) \ \Big| \ 
\ \vc{v}(0, \cdot) = \vc{v}_0, \ \vc{v} \cdot \vc{n}|_{\partial \Omega} = 0, \right. \\ 
&\ \ \Div \vc{v} = 0, \ \partial_t \vc{v} + \Div \mathbb{U} = 0 \ \mbox{for some}\ 
\mathbb{U} \in C^1([0,T] \times \Ov{\Omega}; \R^{d \times d}_{0, {\rm sym}})\\&\left. 
\ \frac{d}{2} \sup_{[0,T] \times \Ov{\Omega}}  \lambda_{\rm max} \left[ \frac{(\vc{v} - h'(t) \Grad \Phi_0 ) \otimes 
(\vc{v} - h'(t) \Grad \Phi_0 )}{\vr} 
- \mathbb{U} \right] < \inf_{[0,T] \times \Ov{\Omega}} e \right\}.
\end{split}
\]
The functions $E$ and $\vc{m}_0$ given, we fix $\Lambda_0$, together with the constant $E_0$, so that the set 
$X_0$ is non--empty. This can be achieved by considering $\vc{v} = \vc{v}_0$, $\mathbb{U} = 0$ and fixing $\Lambda_0$ appropriately. 
Finally, we set 
\[
\Ov{e} = \sup_{[0,T] \times \Ov{\Omega}} e(t,x) < \infty.
\]

Thus, by virtue of \eqref{O1}, the set $X_0$ is bounded in $L^\infty((0,T) \times \Omega; \R^d)$; whence metrizable in the topology 
of $C_{{\rm weak}}([0,T]; L^2(\Omega; \R^d))$. We denote by $X$ its closure in the corresponding metric $d$.

\subsection{Critical points of the energy functional}

Following \cite{DelSze3}, we introduce the functional
\[
I[\vc{v}] = \int_0^T \intO{ \left( \frac{1}{2} \frac{|\vc{v} - h'(t) \Grad \Phi_0|^2}{\vr} - e \right) } \dt \ \mbox{for}\ \vc{v} \in X. 
\]
The functional $I$ is convex lower--semicontinuous on the complete metric space $X$. By Baire category argument we conclude that the points of continuity must form a dense set in $X$.

The second observation is that 
\[
I[\vc{v}] = 0 \ \Rightarrow \ \vc{v} \ \mbox{is a weak solution of the problem \eqref{S15}, \eqref{S16}.}
\]
Indeed, from convexity of the function 
\[
[\vc{v}; \mathbb{U}] \mapsto
\frac{d}{2} \lambda_{\rm max} \left[ \frac{(\vc{v} - h'(t) \Grad \Phi_0) \otimes (\vc{v} - h'(t)\Grad \Phi_0) }{\vr} - \mathbb{U} \right], 
\]
we deduce that for any $\vc{v} \in X$ there is $\mathbb{U} \in L^\infty((0,T) \times \Omega; \R^{d \times d}_{0,{\rm sym}})$
\[
\begin{split}
\partial_t \vc{v} + \Div \mathbb{U} &= 0 \ \mbox{in}\ \mathcal{D}'((0,T) \times \Omega),\\ 
\frac{1}{2} \frac{|\vc{v} - h'(t) \Grad \Phi_0 |^2}{\vr} &\leq \frac{d}{2} \lambda_{\rm max} \left[ \frac{(\vc{v} - h'(t) \Grad \Phi_0) \otimes (\vc{v} - h'(t) \Grad \Phi_0) }{\vr} - \mathbb{U} \right]
\leq e \ \mbox{a.e. in}\ (0,T) \times \Omega.
\end{split}
\]
Consequently, $I \leq 0$ on $X$; while $I[\vc{v}] = 0$ implies the desired relations (cf. \cite{DelSze3})
\[
\begin{split}
\frac{1}{2} &\frac{|\vc{v} - h'(t) \Grad \Phi_0|^2}{\vr} = e,\\ 
\mathbb{U} &= \frac{(\vc{v} - h'(t) \Grad \Phi_0)  \otimes (\vc{v} - h'(t) \Grad \Phi_0)}{\vr} - \frac{d}{2} \frac{|\vc{v} - h'(t) \Grad \Phi_0|^2}{\vr} \mathbb{I} 
\ \mbox{ a.e. in}\ (0,T) \times \Omega.
\end{split}
\] 

Thus, similarly to the arguments used in \cite{DelSze3}, it remains to observe:
\begin{equation} \label{S18}
\vc{v} \ \mbox{-- a point of continuity of}\ I \ \mbox{on}\ X \ \Rightarrow 
I[\vc{v}] = 0.
\end{equation}

To show \eqref{S18}, we argue by contradiction. Assuming 
\[
I[\vc{v}] = \underline{I} < 0
\]
we construct a sequence of functions 
\[
\vc{v}_m \in X_0 \ \mbox{with the corresponding fields}\ \mathbb{U}_m \in C^1([0,T] \times \Ov{\Omega}; \R^{d\times d}_{0, {\rm sym}})
\]
such that
\[
\vc{v}_m \to \vc{v} \ \mbox{in}\ X, \ I[\vc{v}_m] \to \underline{I} < 0\ \mbox{as}\ m \to \infty.
\]

For fixed $m$, we apply Oscillatory Lemma (Lemma \ref{OL3}) for $\vc{v} = \vc{v}_m - h'(t) \Grad \Phi_0$, $\mathbb{U} = \mathbb{U}_m$, 
$r = \vr_0$, and $e$. We obtain sequences $\{ \vc{w}_{m,n} \}_{n=1}^\infty$, $\{ \mathbb{V}_{m,n} \}_{n=1}^\infty$ 
satisfying:
\begin{itemize}
\item
\[
\vc{v}_m + \vc{w}_{m,n} \in X_0 \ \mbox{with the associated fields}\ \mathbb{U}_m + \mathbb{V}_{m,n} 
\ \mbox{for any}\ m,n; 
\]  
\item
\begin{equation} \label{S19}
\vc{v}_m + \vc{w}_{m,n} \to \vc{v}_m \ \mbox{in}\ X \ \mbox{as}\ n \to \infty 
\ \mbox{for any fixed}\ m;
\end{equation}
\item
\begin{equation} \label{S20}
\begin{split}
\liminf_{n \to \infty} \int_0^T\! \!&\intO{ \frac{1}{2} \frac{|\vc{v}_{m} + \vc{w}_{m,n} - h'(t) \Grad \Phi_0 |^2}{\vr} }\dt \\ &= \!\!
\int_0^T \!\!\intO{ \frac{1}{2} \frac{|\vc{v}_m - h'(t) \Grad \Phi_0|^2}{\vr} } + 
\liminf_{n \to \infty} \int_0^T\!\! \intO{ \frac{1}{2} \frac{|\vc{w}_{m,n}|^2}{\vr} }\dt\\
&\geq \int_0^T \intO{ \frac{1}{2} \frac{|\vc{v}_m - h'(t) \Grad \Phi_0 |^2}{\vr} } + c(d, \Ov{e}) 
\int_0^T \intO{ \left( \frac{1}{2} \frac{|\vc{v}_m - h'(t) \Grad \Phi_0|^2}{\vr} - e \right)^2 } \dt\\
&\geq \int_0^T \intO{ \frac{1}{2} \frac{|\vc{v}_m - h'(t) \Grad \Phi_0|^2}{\vr} } + 
c(d, \Ov{e}) T^{-1} |\Omega|^{-1} (I[\vc{v}_m])^2;
\end{split}
\end{equation}
\end{itemize}
where we have used Jensen's inequality in \eqref{S20}. Relation \eqref{S20} rewritten as 
\[
\liminf_{n \to \infty} I [ \vc{v}_m + \vc{w}_{m,n}] \geq I[\vc{v}_m] + 
\frac{c(d, \Ov{e})}{T |\Omega|} \left( I[\vc{v}_m] \right)^2 \ \mbox{for any}\ m
\]
implies that $\vc{v}$ cannot be a point of continuity of $I$ unless $I[\vc{v}] = 0$.

We have proved Theorem \ref{MT1}.

\subsection{Points of strong continuity}

We show how Theorem \ref{MT2} follows from Theorem \ref{MT1}. Given the set $\{ \tau_n \}_{n=1}^\infty$ it is a routine matter to construct a function $E: [0,T] \to \infty$, 
\[
\begin{split}
0 \leq E(t) &\leq \Ov{E} \ \mbox{for all}\ t \in [0,T],\ 
E \ \mbox{strictly decreasing in}\ [0,T],\\ 
&\lim_{t \to \tau_n -} E(t) > \lim_{t \to \tau_n +} E(t) \ \mbox{for any}\ \tau_n, \ n=1,2,\dots
\end{split}
\]

Consider the solutions $[\vr, \vm]$, the existence of which is guaranteed by Theorem \ref{MT1} with the energy profile 
\[
\intO{ \left( \frac{1}{2} \frac{|\vm|^2}{\vr} + P(\vr) \right)(\tau, \cdot) } = E_0 + E(\tau) 
\ \mbox{for a.a.}\ \tau \in (0,T).
\]
As $\vr$, $\vm$ is uniformly bounded and $\vr$ bounded below away from zero, the energy
\[  
\tau \mapsto \intO{ \left( \frac{1}{2} \frac{|\vm|^2}{\vr} + P(\vr) \right)(\tau, \cdot) }
\]
must be continuous at any point of strong continuity of $[\vr, \vm]$. Consequently, $\tau_n$ cannot be points of strong continuity 
of $[\vr, \vm]$. 

We have shown Theorem \ref{MT2}.

\def\cprime{$'$} \def\ocirc#1{\ifmmode\setbox0=\hbox{$#1$}\dimen0=\ht0
  \advance\dimen0 by1pt\rlap{\hbox to\wd0{\hss\raise\dimen0
  \hbox{\hskip.2em$\scriptscriptstyle\circ$}\hss}}#1\else {\accent"17 #1}\fi}


\end{document}